\documentclass[a4paper,12pt,leqno]{article}
\usepackage{latexsym}
\usepackage[all]{xy}

\usepackage{amssymb} 
\usepackage{amsmath} 
\usepackage{theorem}

\def\Z{{\mathbb{Z}}}
\def\K{{\mathbb{K}}}
\def\R{{\mathbb{R}}}
\def\A{{\mathcal{A}}}
\def\B{{\mathcal{B}}}

\DeclareMathOperator{\codim}{codim}

\DeclareMathOperator{\Der}{Der}

\DeclareMathOperator{\Shi}{Shi}
\DeclareMathOperator{\Cat}{Cat}

\numberwithin{equation}{section}

\newcommand{\owari}{\hfill$\square$}

\theoremstyle{break}
\newtheorem{theorem}{Theorem}[section]
\newtheorem{prop}[theorem]{Proposition}
\newtheorem{cor}[theorem]{Corollary}
\newtheorem{lemma}[theorem]{Lemma}
\newtheorem{define}[theorem]{Definition}
\newtheorem{rem}[theorem]{Remark}

\title{Multiple addition, deletion and restriction theorems for 
hyperplane arrangements}
\author{
Takuro Abe\footnote{
Institute of Mathematics for Industry, Kyushu University. abe@imi.kyushu-u.ac.jp} 
and 
Hiroaki Terao
\footnote{
Department of Mathematics, Hokkaido University. terao@math.sci.hokudai.ac.jp}
}
\date{\today} 

\begin{document}

\maketitle

\begin{abstract}
In the study of free arrangements, the most useful result to construct/check 
free arrangements is the addition-deletion theorem in \cite{T}. Recently, the multiple version of 
the addition theorem is proved in \cite{ABCHT}, called the multiple addition theorem (MAT) to prove the ideal-free theorem. The aim of this article is to give the deletion version of MAT, the multiple deletion theorem (MDT). Also, we can generalize MAT from the viewpoint of our new proof. Moreover, we introduce their restriction version, a multiple restriction theorem (MRT). Applications of them including the combinatorial freeness of the 
extended Catalan arrangements are given. 
\end{abstract}

\section{Introduction}
Let $\K$ be an arbitrary field, 
$V=\K^\ell$ and $S=\K[x_1,\ldots,x_\ell]$ its coordinate ring. Let $\Der S:=
\oplus_{i=1}^\ell S \partial_{x_i}$. A \textbf{hyperplane arrangement} $\A$ is a 
finite set of linear hyperplanes in $V$. For each $H \in \A$, let us fix a linear form 
$\alpha_H \in V^*$ such that $\ker \alpha_H=H$. The \textbf{logarithmic derivation module} $D(\A)$ of $\A$ is defined by 
$$
D(\A):=\{\theta \in \Der S \mid 
\theta(\alpha_H) \in S\alpha_H\ (\forall H \in \A)\}.
$$
$D(\A)$ is a reflexive $S$-module, and not free in general. We say that $\A$ is \textbf{free} 
with \textbf{exponents} $\exp(\A)=(d_1,\ldots,d_\ell)$ if there are $S$-free homogeneous basis $\theta_1,
\ldots,\theta_\ell \in D(\A)$ with $\deg \theta_i=d_i\ (i=1,\ldots,\ell)$. For 
integers $d_1,\ldots,d_\ell$, let $(d_1,\ldots,d_\ell)_\le$ denote 
that 
$d_1 \le \cdots \le d_\ell$. 

In the study of hyperplane arrangements,
the freeness is one of the most important objects. But in general it is very difficult to 
determine whether a given arrangment is free or not. The most useful method to 
check the freeness is the following addition-deletion theorem.

\begin{theorem}[\cite{T}]
Let $H \in \A,\ \A':=\A \setminus \{H\}$ and 
$\A'':=\A^H:=\{H \cap L \mid L \in \A'\}$. Then two of the following three imply the third:
\begin{itemize}
\item[(1)]
$\A$ is free with $\exp(\A)=(d_1,\ldots,d_{\ell-1},d_\ell)$.
\item[(2)]
$\A'$ is free with $\exp(\A')=(d_1,\ldots,d_{\ell-1},d_\ell-1)$.
\item[(3)]
$\A''$ is free with $\exp(\A'')=(d_1,\ldots,d_{\ell-1})$.
\end{itemize}
In particular, all the three hold true if $\A$ and $\A'$ are both free.
\label{additiondeletion}
\end{theorem}

There are a lot of variants of Theorem \ref{additiondeletion}. For example, the division theorem in \cite{A} is one of them. Among them, recently, 
the multiple addition theorem (MAT) is introduced in \cite{ABCHT} to prove the freeness of 
ideal subarrangements. Let us recall MAT.

\begin{theorem}[Multiple addition theorem (MAT), \cite{ABCHT}, Theorem 3.1]
Assume that $\A'$ is a free arrangement with 
$\exp(\A')=(d_1,d_2.\ldots,d_\ell)_\le$ such that 
$d_i
=:d$ for $q \le i\le \ell$. For $H_\ell,\ldots,H_q \in \A$ with $2 \le q$, 
define $\A_j'':=(\A' \cup \{H_j\})^{H_j},\ \A:=\A' \cup \{H_q,\ldots,H_\ell\}$ 
and assume the following conditions:
\begin{itemize}
\item[(1)]
$|\A'| -|\A_j''|=d\ (j=q,\ldots,\ell)$, 
\item[(2)]
$X:=\bigcap_{i=q}^\ell H_i$ is $(\ell-q+1)$-codimensional, and 
\item[(3)]
$X \not \subset \bigcup_{K \in \A'}K$. 
\end{itemize}
Then $\A$ is free with exponents 
$(d_1,d_2,\ldots,d_{q-1},d+1,\ldots,d+1)$.
\label{MAT}
\end{theorem}

This multiple addition theorem 
enables us, under certain conditions, to add hyperplanes at once to free arrangements keeping the freeness with controlled behavior of exponents. Then it is very natural to ask, what about the multiple deletion theorem? The first aim of this article is to answer this 
very natural question as follows:

\begin{theorem}[Multiple deletion theorem, MDT]
Assume that $\A$ is a free arrangement with 
$\exp(\A)=(1,d_2.\ldots,d_\ell)_\le$ such that 
$1 < d_2$. For $H_2,\ldots,H_q \in \A$ with $2 \le q \le \ell$, let 
$\alpha_j:=\alpha_{H_j}\ (j=2,\ldots,q)$, 
define $\A_j'':=\A^{H_j},\ \A':=\A \setminus \{H_2,\ldots,H_q\}$ 
and assume the following conditions:
\begin{itemize}
\item[(1)]
$|\A| -|\A_j''|=d_j\ (j=2,\ldots,q)$, 
\item[(2)]
$X:=\bigcap_{i=2}^q H_i$ is $(q-1)$-codimensional, and 
\item[(3)]
$X \not \subset \bigcup_{K \in \A'}K$. 
\end{itemize}
Then $\A'$ is free with exponents 
$(1,d_2-1,d_3-1,\ldots,d_q-1,d_{q+1},\ldots,d_\ell)$. More explicitly, 
there is a basis for $D(\A)$ of the form $\theta_E,\alpha_2 \eta_2,\ldots,
\alpha_q \eta_q,\theta_{q+1},\ldots,\theta_\ell$ for $D(\A)$ with 
$\deg \eta_i=d_i-1,\ \deg \theta_j=d_j$ for all $i,j$.
\label{MDT2}
\end{theorem}

From the same viewpoint of Theorem \ref{MDT2}, we can show a generalized version of 
the multiple addition theorem in \cite{ABCHT} as follows:

\begin{theorem}[Multiple addition theorem 2, MAT2]
Assume that $\A'$ is a free arrangement with 
$\exp(\A')=(1,d_2.\ldots,d_\ell)_\le$. For $H_\ell,\ldots,H_q \in \A$ with $2 \le q$, 
define $\A_j'':=(\A' \cup \{H_j\})^{H_j},\ \A:=\A' \cup \{H_q,\ldots,H_\ell\}$ 
and assume the following conditions:
\begin{itemize}
\item[(1)]
$|\A'| -|\A_j''|=d_j\ (j=q,\ldots,\ell)$, 
\item[(2)]
$X:=\bigcap_{i=q}^\ell H_i$ is $(\ell-q+1)$-codimensional, and 
\item[(3)]
$X \not \subset \bigcup_{K \in \A'}K$. 
\end{itemize}
Then $\A$ is free with exponents 
$(1,d_2,\ldots,d_{q-1},d_{q}+1,\ldots,d_\ell+1)$. Moreover, there is a basis 
$\theta_E,\theta_2,\ldots,\theta_{q-1},\eta_q,\ldots,\eta_\ell$ for $D(\A')$ such that 
$\deg \theta_i=d_i,\ \deg \eta_j=d_j$, 
$\theta_i \in D(\A)$ and $\eta_j \in D(\A \setminus \{H_j\})$ for all $i,\ j$.
\label{MAT2}
\end{theorem}

Also, we can show the restriction version as follows:

\begin{theorem}[Multiple restriction theorem, MRT]
Let $\A$ be free with $\exp(\A)=(1,d_2,\ldots,d_\ell)_\le$ and let 
$X \in L_k(\A)$. Then $\A^X:=\{X \cap H \mid H \in \A,\ X \not \subset H\}$ is free with $\exp(\A^X)=(1,d_{k+2},\ldots,d_\ell)$ 
if and only if $|\A^X|=1+d_{k+2}+\cdots+d_{\ell}$.
\label{cor1}
\end{theorem}

We give some applications of MDT and MAT2. In particular, by observing MAT, we can show that 
the freeness of the extended Catalan arrangements depends only on its combinatorics. 

The organization of this article is as follows. In \S2 we introduce several definitions and results used in this article. In \S3 we prove MDT and MAT2. In \S4 we prove MRT by 
showing a more general result. In \S5 we consider the relation between combinatorics of freeness and MAT.
\medskip

\noindent
\textbf{Acknowledgements}. 
The first author is 
partially supported by KAKENHI, JSPS Grant-in-Aid for Scientific Research (B) 16H03924, and 
Grant-in-Aid for Exploratory Research 16K13744.

\section{Preliminaries}
In this section let us fix a notation and introduce several results used for the rest of this 
article. Let $\A$ be an arrangement of hyperplanes in $V=\K^\ell$ as in \S1. 
The \textbf{intersection lattice} $L(\A)$ of $\A$ is defined by 
$$
L(\A):=\{\bigcap_{H \in \B} H \mid 
\B \subset \A\}.
$$
Let $L_k(\A):=\{ X \in L(\A) \mid \codim_V(X)=k\}$. 
On $L(\A)$ the \textbf{M\"obius function} $\mu$ is defined by $\mu(V)=1$, and by 
$$
\mu(X)=-\sum_{Y \in L(\A),\ X \subsetneq Y \subset V} \mu(Y).
$$
The \textbf{Poincar\'e polynomial} $\pi(\A;t)$ is 
$$
\pi(\A;t):=\sum_{X \in L(\A)} \mu(X)(-t)^{\codim X},$$
and the \textbf{characteristic polynomial} $\chi(\A;t)$ is 
$$
\chi(\A;t):=\sum_{X \in L(\A)} \mu(X)t^{\dim X}.
$$
For $X \in L(\A)$, the \textbf{localizaton} $\A_X:=\{H \in \A 
\mid X \subset H\}$ is a subarrangement of $\A$, and 
the \textbf{restriction} $\A^X:=\{H \cap X \mid H \in \A \setminus \A_X\}$ 
is an arrangement in $X$. More generally, for a point $x \in V$, let $\A_x:=\{
H \in \A \mid x \in H\}$, and $D(\A)_x$ is the localization of $D(\A)$ at the 
homogeneous prime ideal corresponding to $x$.

Next recall several results related to $D(\A)$ introduced in \S1. 
For a non-negative integer $d$, let $D(\A)_d$ denote the homogeneous degree $d$-part of 
$D(\A)$. Here $\theta \in \Der S$ is homogeneous of degree $d$ if $\theta(\alpha)$ is 
homogeneous of degree $d$ for all $\alpha \in V^*$ with $\theta(\alpha) \neq 0$.
When $\A \neq \emptyset$, $\theta_E \in D(\A)$ and the free $S$-module $S\theta_E$ 
is a direct summand of $D(\A)$. Explicitly, for any $H \in \A$, 
if $D_H(\A):=\{\theta \in D(\A) \mid \theta(\alpha_H)=0\}$, then 
$$
D(\A) \simeq S\theta_E \oplus D_H(\A).
$$
Hence $\A$ is free if and only if $D_H(\A)$ is free for some, and hence any $H \in \A$. 
Hence $1 \in \exp(\A)$ when $\A \neq \emptyset$ is 
free. For $X \in L(\A)$, it is easy to see that $\A_X$ is free if $\A$ is free. Contrary to 
the localization, there are no relations of freeness between $\A$ and $\A^H$ for $H \in \A$. 
The most useful relation is Theorem \ref{additiondeletion}. In fact, always we have the 
exact sequence 
\begin{equation}
0 \rightarrow D(\A') \stackrel{\cdot \alpha_H}{\rightarrow} D(\A) 
\stackrel{\rho}{\rightarrow} D(\A''),
\label{e}
\end{equation}
where $\rho$ is the map taking modulo $\alpha_H$, and called the Euler restriction map. And the exact sequence (\ref{e}) is 
right exact if all the three conditions in Theorem \ref{additiondeletion} hold.
When $\A$ is free, we can relate the exponents with its Poincar\'e polynomial as follows.

\begin{theorem}[\cite{T1}]
If $\A$ is free with $\exp(\A)=(d_1,\ldots,d_\ell)$, then $\pi(\A;t)=\prod_{i=1}^\ell 
(1+d_i t)$.
\label{factorization}
\end{theorem}

For a \textbf{multiplicity function} $m:\A \rightarrow \Z_{>0}$, we can define 
the \textbf{logarithmic derivation module} $D(\A,m)$ of the multiarrangement $(\A,m)$ as 
$$
D(\A,m):=\{\theta \in \Der S \mid \theta(\alpha_H) \in S\alpha_H^{m(H)}\ (\forall H \in \A)\}.
$$
$D(\A,m)$ is a reflexive $S$-module, and not free in general. Hence we can define its 
\textbf{freeness} and 
\textbf{exponents} in the same way as for $\A$. Let $|m|:=\sum_{H \in \A} m(H)$ and 
$Q(\A,m):=\prod_{H \in \A} \alpha_H^{m(H)}$.

\begin{define}
Let $H \in \A$. Then the \textbf{Ziegler restriction} $(\A^H,m^H)$ of $\A$ onto $H$ is 
defined by 
$$
m^H(X):=|\{L \in \A \setminus \{H\} \mid L \cap H=X\}|
$$
for $X \in L(\A^H)$. 
Then we have the \textbf{Ziegler restriction map} 
$$
\pi_H:D_H(\A) \rightarrow D(\A^H,m^H)$$
by taking modulo $\alpha_H$. Hence $\ker \pi_H=D_H(\A)[-1]=\alpha_HD_H(\A)$.
\label{Zrest}
\end{define}

\begin{theorem}[\cite{Z}]
Assume that $\A$ is free with $\exp(\A)=(1,d_2,\ldots,d_\ell)$. Then $(\A^H,m^H)$ is free 
with $\exp(\A^H,m^H)=(d_2,\ldots,d_\ell)$. In particular, the Ziegler restriction map 
is surjective in this case.
\label{Zieglerfree}
\end{theorem}

For the freeness, the following criterion is important.

\begin{theorem}[Saito's criterion, \cite{S}]
Let $\theta_1,\ldots,\theta_\ell \in D(\A,m)$ be homogeneous 
derivations. Then 
$\det (\theta_i(x_j)) \in S Q(\A,m)$. Moreover, they 
form a 
basis for $D(\A,m)$ if and only if 
they are $S$-independent, and $|m|=\sum_{i=1}^\ell \deg \theta_i$.
\label{Saito}
\end{theorem}

The following is very important to analyze freeness.

\begin{theorem}[\cite{T}]
Let $H \in \A$ and $\A':=\A \setminus \{H\}$. Then there is a polynomial $B$ with 
$\deg B=|\A'|-|\A^H|$ such that 
$$
\theta(\alpha_H) \in (\alpha_H,B)
$$
for all $\theta \in D(\A')$.
\label{polyB}
\end{theorem}

Based on Theorem \ref{polyB}, we often use the following two arguments in this article.

\begin{prop}
Let $H \in \A,\ \A':=\A \setminus \{H\}$ and 
$d:=|\A'|-|\A^H|$. Then 

(1)\,\,
$\theta \in D(\A')$ belongs to $D(\A)$ if $\deg \theta <d$.

(2)\,\,
Assume that there is $\theta \in D(\A')_d$ such that 
$\theta(\alpha_H) =B$ modulo $\alpha_H$. 
If there is a generator $\theta_E,\varphi_1,\ldots,\varphi_s, \theta$ for 
$D(\A')$ such that $\varphi_i(\alpha_H)\equiv g_i B$ for $i=1,\ldots,s$ 
modulo $\alpha_H$, then 
$\theta_E,\varphi_1-g_i \theta,\ldots,\varphi_s-g_s \theta, \alpha_H \theta$ is a generator for 
$D(\A)$. In particular, 

(3)\,\,
if $\A'$ is free with $\exp(\A')=(d_1,\ldots,d_\ell),\ d_\ell=d$,
$\theta_1,\ldots,\theta_\ell$ form a basis 
for $D(\A')$ with $\deg \theta_i=d_i$, $\theta_\ell(\alpha_H) =B$ modulo $\alpha_H$,  
and $\theta_i(\alpha_H) \equiv g_iB$ modulo $\alpha_H$, then 
$\theta_1-g_1 \theta,\ldots,\theta_{\ell-1}-g_{\ell-1} \theta,\alpha_H \theta_\ell$ form a basis for $D(\A)$.
\label{Bappli}
\end{prop}

\noindent
\textbf{Proof}.
(1)\,\, Trivial by Theorem \ref{polyB}. 

(2)\,\,
For $\varphi_i$ with $\deg \varphi_i < d$, (1) implies that $\varphi_i \in D(\A)$. For 
those with $\deg \varphi_i \ge d$, let 
$$
\varphi_i(\alpha_H)=f_i \alpha_H + g_i B
$$
by Theorem \ref{polyB}. Then replacing $\varphi_i$ by $\varphi_i-g_i \theta$, 
it holds that all $\varphi_i \in D(\A)$. Let us show that $\theta_E,\varphi_1,\ldots,
\varphi_s,\alpha_H\theta$ generate $D(\A)$. Note that 
$\theta(\alpha_H) \not \in S \alpha_H$. Let $\varphi \in D(\A) \subset D(\A')$. Then 
$\varphi=f \theta_E+\sum_{i=1}^s f_i \varphi_i+g\theta$. 
Since all derivations but $\theta$ belongs to $D(\A)$, $g \theta \in D(\A)$ which implies that 
$\alpha_H \mid g$. Thus (2) follows. 
(3) is immediate from (2).\owari
\medskip

\begin{prop}
Let $H \in \A$ and assume that $\A$ is free with $\exp(\A)=(1,d_2,\ldots,d_\ell)$. Let 
$\theta_i\ (i=2,\ldots,\ell)$ be a basis for $D(\A)$ with $\deg \theta_i=d_i$, and let 
$Q':=Q(\A^H,m^H)/Q(\A^H)$.  
Assume that $|\A|-|\A^H|=d_\ell$, and 
$$
Q'\pi_H(\theta_E)=\sum_{i=2}^\ell
\pi_H(f_i) \pi_H (\theta_i),
$$
where $f_i \in S$. If there is some $i$ such that $\pi_H(f_i) \neq 0$ and 
$d_i=d_\ell$, then 
$$
\pi_H(\theta_2),\ldots,\pi_H(\theta_{\ell-1}),Q' \pi_H(\theta_E)
$$
form a basis for $D(\A^H,m^H)$, and 
$$
\pi(\theta_E),\pi_H(\theta_2),\ldots,\pi_H(\theta_{\ell-1})
$$
form a basis for $D(\A^H)$.
\label{replace}
\end{prop}

\noindent
\textbf{Proof}. 
By Theorem \ref{Zieglerfree} and the decomposition $D(\A) \simeq S\theta_E 
\oplus D_H(\A)$, we may assume that  
$\pi_H(\theta_2),\ldots,\pi_H(\theta_{\ell-1}),\pi_H(\theta_\ell)$ form a basis for $D(\A^H,m^H)$. We may assume that $f_\ell \neq 0$. By the reason of degrees, 
$f_\ell \in \K^\times$. Hence we may replace $\pi_H(\theta_\ell)$ by $Q'\pi_H(\theta_E)$ to 
obtain that 
$\pi_H(\theta_2),\ldots,\pi_H(\theta_{\ell-1}),Q' \pi_H(\theta_E)$ is a basis for $D(\A^H,m^H)$. In particular, 
$\pi_H(\theta_E), \pi_H(\theta_2),\ldots,\pi_H(\theta_{\ell-1})$ are independent over $S/\alpha_H$, and the sum of their degrees are $|\A^H|$. Hence Saito's criterion completes the proof.\owari
\medskip

\section{Proof of Theorems \ref{MDT2} and \ref{MAT2}}

Let us prove main results in this 
article.
\medskip

\noindent
\textbf{Proof of Theorem \ref{MDT2}}.
Let $\alpha_i$ be the defining linear form of $H_i$. We show by induction on $q$. First let $q=2$.  Let $\theta_E,\theta_2,\ldots,\theta_\ell$ be a basis for $D(\A)$ with $\deg \theta_i=d_i$ and 
$\theta_i \in D_{H_2}(\A)$. Let $\pi:D_{H_2}(\A) \rightarrow D(\A''_2,m^{H_2})$ be the Ziegler restriction map. Then Theorem \ref{Zieglerfree} implies that 
$\pi(\theta_2),\ldots,\pi(\theta_\ell)$ form a basis for $D(\A''_2,m^{H_2})$. Let $Q:=Q(\A''_2,m^{H_2})/Q(\A''_2)$. Then $Q\pi(\theta_E) \in D(\A''_2,m^{H_2})_{d_2}$ by the condition (1). Hence 
there are $f_i \in S$ such that 
$$
Q\pi(\theta_E)=\sum_{i=2}^k \pi(f_i) \pi(\theta_i),
$$
where $d_2=\cdots=d_k<d_{k+1}$. 
By comparing the degrees, we may assume that 
$\pi(f_2) \in \K^\times$. 
Then by Proposition \ref{replace}, we may replace $\pi(\theta_2)$ by $Q\pi(\theta_E)$ to obtain a basis 
$Q\pi(\theta_E),\pi(\theta_3),\ldots,\pi(\theta_\ell)$ for $D(\A''_2,m^{H_2})$. Hence again Proposition \ref{replace} implies that 
the derivations $\pi(\theta_E),\pi(\theta_3),\ldots,\pi(\theta_\ell)$ form a basis for $D(\A_2'')$. 
Let $\rho:D(\A) \rightarrow D(\A_2'')$ be the Euler restriction map. Then $\rho(\theta_i)=
\pi(\theta_i)$ for $i \ge 3$ by definition. Since $\rho$ is surjective too, 
there are polynomials $f,\ f_i\ (i=3,\ldots,\ell)$ such that $$
\theta_2=f \theta_E+\sum_{i=3}^\ell f_i\theta_i+\alpha_H \theta'
$$
by (\ref{e}). Here $\theta' \in D(\A \setminus \{H_2\})$. Replacing 
$\theta_2$ by $\theta_2-f\theta_E-\sum_{i=3}^\ell f_i \theta_i$, we have the desired basis 
for $D(\A)$. 

Now assume that the statement holds true for 
$q-1 \ge 1$. When the case of $q$, by induction, $D(\A)$ has a basis of the form 
$\theta_E,\alpha_2 \eta_2,\ldots,\alpha_{q-1} \eta_{q-1}, \theta_q,\ldots,\theta_\ell$ with 
$\deg \eta_j=d_j-1,\ \deg \theta_j=d_j$.
Replacing $\alpha_j \eta_j$ by $\alpha_j \eta_j - (\alpha_j \eta_j(\alpha_q)/\alpha_q)\theta_E$, we may assume that $\alpha_j \eta_j \in D_{H_q}(\A)$ for $j=2,\ldots,q-1$. Also, we may assume that $\theta_q,\ldots,\theta_\ell \in D_{H_q}(\A)$. 
Let $\pi:D_{H_q}(\A) \rightarrow D(\A^{H_q},m^{H_q})$ be 
the Ziegler restriction map.
By Theorem \ref{Zieglerfree}, $\pi(\alpha_2 \eta_2),\ldots,\pi(\alpha_{q-1} 
\eta_{q-1}),\pi(\theta_q),\ldots,\pi(\theta_\ell)$ form a basis for $D(\A^{H_q},m^{H_q})$.
Let $Q:=Q(\A^{H_q},m^{H_q})/Q(\A^{H_q})$. Then $\theta_E^q:=
Q\pi(\theta_E) \in D(\A^{H_q},m^{H_q})_{d_q}$ by the condition 
(1). Hence $\theta_E^q$ can be expressed as a $S/\alpha_q$-linear combination of the basis whose degree is at most $d_q$. If there is $\theta_j$ such that $\deg \theta_j=d_q$, and it appears in the linear combination with a non-zero scalar coefficient, then the same argument 
as when 
$q=2$, combined with 
Proposition \ref{replace}, implies that $\pi(\alpha_2 \eta_2),\ldots,\pi(\alpha_{q-1} 
\eta_{q-1}),\theta_E^q,\pi(\theta_{q+1}),
\ldots,\pi(\theta_\ell)$ 
form a basis for $D(\A''_q,m^{H_q})$. Hence 
by the same argument as when $q=2$, we may replace $\theta_q$ by the derivation 
of the form $\alpha_{H_q} \theta'$ with $\theta' \in D(\A \setminus \{H_q\})$, which 
completes the proof. 
Assume not. Hence 
\begin{equation}\label{eq100}
\theta_E^q=\sum_{i=2}^{q-1} \pi(f_i) \pi(\alpha_i \eta_i)\ (f_i \in S).
\end{equation}
Let $H_q \cap H_j=:X_j \in \A^{H_q}$ and let $\beta_j \in S/\alpha_q$ be the defining linear form 
of $X_j$. Note that $\cap_{j=2}^{q-1} X_j =X$. 
We show that there is a point $x \in X$ such that $\alpha_Y(x) \neq 0$ 
for all $Y \in \A^{H_q} \setminus \{X_2,\ldots,X_{q-1}\}$. To show it, it suffices to check that 
$X \not \subset Y\in \A^{H_q} \setminus \{X_2,\ldots,X_{q-1}\}$ by extending the base 
field $\K$ if necessary.
Assume that there is such $Y=H \cap H_q,\ H \in \A$. 
Since $X \subset Y \subset H$, the condition (3) implies that 
$H \in \{H_2,\ldots,H_{q-1}\}$, a contradiction. 

Hence there is $x \in X \setminus \bigcup_{Y \in \A''_q,\ Y \neq X_j} Y$ such that 
the right hand side of (\ref{eq100}) is zero at $x$. We show that 
the left hand side cannot be zero at $x$. 
Since $\pi(\theta_E)$ is nowhere vanishing, it is equivalent to say that 
$Q(x)=0$. 
Note that $\alpha_Y(x) \neq 0$ for $Y \in 
\A''_q \setminus \{X_2,\ldots,X_{q-1}\}$. Thus $Q(x)=0$ could occur only when 
$m(X_j) \ge 2$ for some $j$, which cannot occur by the 
condition (2). 
%
\medskip

\noindent
\textbf{Proof of Theorem \ref{MAT2}}.
Let $\alpha_i$ be the defining linear form of $H_i$. We show by induction on $q$. When 
$q=\ell$,
this is Theorem \ref{MAT}. 
Assume that the statement holds true for 
$q+1\le \ell$. By induction, $D(\A')$ has a basis of the form 
$\theta_E,\theta_2,\ldots,\theta_q,\eta_{q+1},\ldots,\eta_\ell$ with 
$\deg \theta_i=d_i,\ \deg \eta_j=d_j$ such that 
$\eta_j(\alpha_j)=b_j$ modulo $\alpha_j$ for $j=q+1,\ldots,\ell$, where 
$b_j$ is the polynomial in Theorem \ref{polyB} with $\deg b_j=|\A'|-|\A''_j|$. 
Also, if $d_i < d_q$, then Proposition \ref{Bappli} implies that 
$\theta_i \in D(\A)$. 
Moreover, by the induction hypothesis, $\eta_j(\alpha_i)\equiv 0$ modulo $\alpha_i$ for all $ i \neq j$. 
Hence it suffices to show that $\theta_i(\alpha_q) \equiv b_q$ modulo $\alpha_q$ for some $i$ with 
$\deg \theta_i=d_q$ by Proposition \ref{Bappli}. Assume not. 
Then we may assume that $\theta_i$ is tangent to $X$ 
for all $i$. Let $x \in X$ with $\A_x=\{H_q,\ldots,H_\ell\}$ by the condition (3). 
Then $$
k(x) \otimes D(\A')_x \simeq \K^\ell =T_{X,x} \oplus N_{X,x},
$$
i.e., the tangent space decomposes into the tangent space and normal 
space of $X$ at $x$. Here $k(x)$ is the residue field of $S_x$. By the condition (2), the former is of $(q-1)$-dimensional 
and the latter $(\ell-q+1)$-dimensional. Note that $\theta_E$ and  
all $\theta_i \in T_{X,x}$ by the 
assumption above. Thus 
\begin{eqnarray*}
\ell-q+1&=&\dim N_{X,x}\\
&=&\dim \langle \eta_{q+1},\ldots,\eta_\ell\rangle_\K\\
&\le&\ell-q,
\end{eqnarray*}
a contradiction.\owari
\medskip

In fact, the assertions in Theorems \ref{MDT2} and \ref{MAT2} are weaker than 
proved here. Explicitly, the following holds.

\begin{cor}
Let $\A$ be a free arrangement satisfying the conditions (1), (2) and (3) in Theorem \ref{MDT2}. 
Let $I \subset \{2,\ldots,q\}$ and set 
$\A_{-I}:=\A \setminus \{H_i\}_{i \in I}$. Then $\A_{-I}$ is free with basis 
$$
\theta_E,\{\eta_i\}_{i \in I},\ \{\alpha_i \eta_i\}_{i \not \in I},\ 
\theta_{q+1},
\ldots,
\theta_\ell
$$
in the notation of Theorem \ref{MDT2}.
\label{MDTexplicit}
\end{cor}

\noindent
\textbf{Proof}. Immediately from the proofs of Theorem \ref{MDT2} and 
Theorem \ref{Saito}. \owari
\medskip

\begin{cor}
Let $\A'$ be a free arrangement satisfying the conditions (1), (2) and (3) in Theorem \ref{MAT2}. 
Let $I \subset \{q,\ldots,\ell\}$ and set 
$\A_{+I}:=\A \cup \{H_i\}_{i \in I}$. Then $\A_{+I}$ is free with basis 
$$
\theta_E,\{\eta_i\}_{i \not \in I},\ \{\alpha_i \eta_i\}_{i  \in I},\ 
\theta_{2},
\ldots,
\theta_{q-1}
$$
in the notation of Theorem \ref{MAT2}.
\label{MATexplicit}
\end{cor}

\noindent
\textbf{Proof}. Immediate by the proofs of Theorem \ref{MAT2} and 
Theorem \ref{Saito}. \owari
\medskip

\section{Multiple restriction theorems (MRT)}

After proving Theorems \ref{MDT2} and \ref{MAT2}, it is also natural to ask whether 
we have a multiple restriction theorem or not. In the setup of MDT and MAT2, we can say yes 
to this question.

For an $\ell$-arrangement $\A$ and a subspace $X \subset V$, define the arrangement 
$\A \cap X$ in $X$ by 
$$
\A \cap X:=\{L \cap X \mid L \in \A,\ L \not \supset X\}.
$$
If $ X \in L(\A)$, then $\A \cap X=\A^X$. By using these terminology, first, let us 
formulate the multiple restriction theorem in terms of MAT2 and MDT. 

\begin{theorem}[MRT from MDT]
Let $\A$ be a free arrangement satisfying the conditions (1), (2) and (3) in 
Theorem \ref{MDT2}. For $I \subset \{2,\ldots,q\}$, let 
$X_I:=\cap_{i \in I} H_i$. 
Then $\A^{X_I}$ is free with basis 
$$
\overline{\theta}_E,
\{\overline{\alpha_i \eta}_i\}_{ i \not \in I},\ 
\overline{\theta}_{q+1},
\ldots,\overline{\theta}_\ell.
$$
Here for $\theta \in D(\A)$, $\overline{\theta}$ stands for the restriction of $\theta$ onto $X$.
\label{MRT1}
\end{theorem}

\noindent
\textbf{Proof}. Apply Corollary \ref{MDTexplicit}. \owari
\medskip

\begin{theorem}[MRT from MAT2]
Let $\A'$ be a free arrangement satisfying the conditions (1), (2) and (3) in 
Theorem \ref{MAT2}. For $I \subset \{q,\ldots,\ell\}$, let 
$X_I:=\cap_{i \in I} H_i$. 
Then $\A' \cap X_I$ is free with basis 
$$
\overline{\theta}_E,
\{\overline{\eta}_i\}_{ i \not \in I},\ 
\overline{\theta}_{2},
\ldots,\overline{\theta}_q.
$$
Here for $\theta \in D(\A')$, $\overline{\theta}$ stands for the restriction of $\theta$ onto $X$.
\label{MRT2}
\end{theorem}

\noindent
\textbf{Proof}. Apply Corollary \ref{MATexplicit}. \owari
\medskip

In fact, when we want to check the freeness of $\A^X$ in terms of Theorem \ref{MRT1}, 
sometimes the conditions 
(1), (2) and (3) are too many to check. We can show the following MRT.




\begin{theorem}
Let $\A$ be an $\ell$-arrangement with 
a generator 
$\theta_E,\theta_1,\ldots,\theta_s,\varphi_1,\ldots,\varphi_{\ell-k-1}$ such that
$\deg \theta_i=d_i\le d_{i+1} \le e_j:=\deg \varphi_j$ for all $i,\ j$. Let $X \in L_k(\A)$ such that 
$|\A^X|=1+\sum_{i=1}^{\ell-k-1} e_i$. Then $\A^X$ is free with 
$\exp (\A^X)= (1,e_1,\ldots,e_{\ell-k-1})$.
\label{rest2}
\end{theorem}

To prove Theorem \ref{rest2}, the fundamental lemma is the following.

\begin{lemma}
Let $\A$ be an arrangement and $X \in L_k(\A)$. Then the Euler 
restriction map $\rho:D(\A) \rightarrow D(\A^X)$ is generically surjective. 
In particular, if $\theta_E,\theta_1,\ldots,\theta_s$ is a generator for $D(\A)$, then 
there are $(\ell-k-1)$-elements among $\theta_1,\ldots,\theta_{\ell-k-1}$ which are 
$S^X$-independent, where $S^X:=\mbox{Sym}^*(X^*)$.
\label{indep}
\end{lemma}

\noindent
\textbf{Proof}. 
Let $x \in X$ be a generic point. Then it is clear that 
\begin{eqnarray*}
D(\A)_x&=&\langle \theta_1,\ldots,\theta_k,
\partial_{x_{k+1}},\ldots,\partial_{x_\ell}\rangle,\\
D(\A^X)_x&=&\langle \partial_{x_{k+1}},\ldots, \partial_{x_\ell}\rangle,
\end{eqnarray*}
where we change the coordinates $\{x_i\}_{i=1}^\ell$ in such a way that $X=\{x_1=\cdots=x_{k}=0\}$ and $\theta_i 
\in D(\A_X) \cap (\oplus_{i=1}^k S \partial_{x_i})$. Hence $\theta_i|_X=0$ for all $i$. 
At this point 
$\rho$ is clearly surjective. Since $\rho$ is an $S$-module homomorphism, 
the image of $\rho$ has to contain at least $(\ell-k)$-elements which contains 
the Euler derivation, and they are $S^X$-independent. 

Also, note that, we may define the map $\rho_{H,Y}^X:D_H(\A) \rightarrow D_Y(\A^X)$ for 
$H \in \A$ with $H \supset X$ and $Y \in \A^X$ as follows. First, for $\theta \in D_H(\A)$ and 
$L \in \A^H$, 
define the map $\rho_{H,L}^H:D_H(\A) \rightarrow D_L(\A^H)$ as 
$$
\rho_{H,L}^H(\theta):=\rho(\theta)-\frac{\rho(\theta)(\alpha_L)}{\alpha_L} \rho(\theta_E),
$$
where 
$\rho$ is the Euler restriction map. Combine this map till $X$ to obtain the 
map $\rho_{H,Y}^X$. Then it is clear that 
$$
\rho=\rho_E\oplus \rho_{H,Y}^X:D(\A)=S \theta_E \oplus D_H(\A) \rightarrow 
S^X \rho_E(\theta_E) \oplus D_Y(\A^X).
$$
Hence the same result holds true for the generator of $D_H(\A)$ and their 
images in $D_Y(\A^X)$.
\owari
\medskip

\noindent
\textbf{Proof of Theorem \ref{rest2}}. 
Assume that the image of $\varphi_1,\ldots,\varphi_{\ell-k-1}$ are $S^X$-independent in 
$D_Y(\A^X)$ in the terminology of the proof of Theorem \ref{rest2}. Then Saito's 
criterion completes the proof. Assume not. Then at least one of $\theta_i$ has to be 
contained in the $S^X$-independent images of the generator. Let 
$M$ be the matrix of coefficients of those $S^X$-independent derivations in 
$D_Y(\A^X)$. Then by Theorem \ref{Saito}, 
$$
|\A^X|-1 \le \deg \det M \le |\A^X|-1-e_j+d_i \le |\A^X|-1
$$
for some $i,\ j$. 
Hence all the inequalities above are equal. Thus again Saito's criterion completes the proof. \owari
\medskip

Now Theorem \ref{cor1} is an immediate corollary of Theorem \ref{rest2}. 



\section{Combinatorial freeness and MAT2}

In this section, we investigate the relation between combinatorial freeness and 
MAT. Though MAT was introduced in \cite{ABCHT} to check freeness of some arrangements, 
it also works well to check whether the freeness is combinatorial or not. Explicitly, 
we can show the following. 

\begin{theorem}
Assume that $\A$ is free and its freeness depends only on $L(\A)$. 
If $\B$ is a free arrangement constructed by adding hyperplanes to $\A$ 
using MAT2, then the freeness of $\B$ depends only on $L(\B)$. 
\label{TC}
\end{theorem}

\noindent
\textbf{Proof}. 
Whether we can add hyperplanes to $\A$, without destroying the freeness, by using MAT2 or not 
keeping freeness depends only on $L(\B)$. This observation completes the proof. \owari
\medskip

To give an example of Theorem \ref{TC}, let us introduce the notation used in this section. 
Let $\K=\R$ and let $\Phi$ an irreducible crystallographic root system of rank $\ell$. Let $W$ be 
the corresponding Weyl group acting on $V=\R^\ell$. 
Fix a positive system $\Phi^+$. Then 
the simple system $\Delta=\{\alpha_1,\ldots,\alpha_\ell\}$ is also fixed. 
For each $\alpha \in \Phi^+$, the hyperplane $H_\alpha$ consists of the points 
$x \in V$ orthogonal to $\alpha$. A \textbf{Weyl arrangement} $\A_{\Phi^+}$ is defined as the set 
$\{H_\alpha \mid \alpha \in \Phi^+\}$. More generally, for $U \subset 
\Phi^+$, define 
$$
\A_U:=\{H_\alpha 
\mid \alpha \in U\}.
$$
Let us introduce a partial order in $\Phi^+$ as follows. 
For $\alpha,\beta \in \Phi^+$, $\alpha \le \beta$ if $
\beta-\alpha \in \sum_{i=1}^\ell \Z_{\ge 0} \alpha_i$. Then a 
\textbf{(lower) ideal} $I \subset \Phi^+$ is the lower closed set with respect to 
this partial order in $\Phi^+$. The main result in \cite{ABCHT}, which was proved by 
using MAT, asserts that $\A_I$ is free if $I$ is an ideal. Also, for $\alpha \in \Phi^+$ and 
$j \in \Z$, define a hyperplane $H_\alpha^j$ in $\R^{\ell+1}=\mbox{Spec}(S[z])$ by 
$$
H_\alpha^j:=\{\alpha=jz\}.
$$
Then the \textbf{extended Shi arrangement} $\Shi^k$ is the arrangement in 
$\R^{\ell+1}$ defined as 
$$
\Shi^k:=
\{H_\alpha^j \mid \alpha \in \Phi^+,\ 
j \in \Z,\ -k+1 \le j \le k\} \cup \{z=0\},
$$
and the \textbf{extended Catalan arrangement} $\Cat^k$ is the arrangement in 
$\R^{\ell+1}$ defined as 
$$
\Cat^k:=
\{H_\alpha^j \mid \alpha \in \Phi^+,\ 
j \in \Z,\ -k \le j \le k\} \cup \{z=0\}.
$$
Both of them are free, which was shown in \cite{Y}.
Now we can give a main application of MAT2. 

\begin{theorem}
The freeness of the extended Catalan arrangements depends only on 
its intersection lattices.
\label{Catalan}
\end{theorem}

Theorem \ref{Catalan} follows immediately from the following.

\begin{cor}
Let $I$ be a lower ideal of the positive system, and 
$$
\Shi^k_{+I}:=
\Shi^k \cup \{H_\alpha^{-k} \mid \alpha \in I\}.
$$
Then the freeness of $\Shi^k_{+I}$ depends only on 
$L(\Shi^k_{+I})$. 
\label{ideal}
\end{cor}

\noindent
\textbf{Proof}. 
We know that the freeness of $\Shi^k$ depends only on the combinatorics by 
Theorem 6.1 in \cite{A}. Hence by Theorem \ref{TC}, it suffices to show that 
$\Shi^k_{+I}$ can be constructed from $\Shi^k$ keeping freeness by using MAT2. 

By Theorem 1.2 in \cite{AT}, we know that 
$\Shi^k_{+I}$ is free with 
$$
\exp_0(\Shi^k_{+I})=
((kh)^\ell)+(d_1,\ldots,d_\ell),
$$
where $(d_1,\ldots,d_\ell)$ is the dual partition of the height 
distribution of roots in $I$ (see \cite{ABCHT}), and 
$\exp_0(\Shi^k_{
+I})$ indicates the multiset of degrees of free basis of 
$D_0(\Shi^k_{
+I}):=D(\Shi^k_{
+I})/S\theta_E$. Let $s$ be the largest height of the positive root belonging to $I$, and 
let 
$J$ be the ideal consisting of positive roots in $I$ whose heights are 
strictly less than 
$s$. We use the induction on the largest height of ideals to show the statement. 
When $s=0$ then there is nothing to show. Assume $s>0$. Then Theorem 1.3 in \cite{AT} implies that 
$\exp_0(\Shi^k_{
+J})=((kh)^\ell)+(d_1,\ldots,d_{\ell-t},(s-1)^t)$, where 
$\{\beta_1,\ldots,\beta_t\}$ is all the positive roots in $I$ of height $s$ and 
$d_1 \le \cdots \le d_{\ell-t} \le s-1$. By Theorems 1.2 and 1.3 in \cite{AT} 
again, $\Shi^k_{
+J} \cup \{H_{\beta_i}^{-k}\}$ is free with exponents 
$((kh)^\ell)+(d_1,\ldots,d_{\ell-t},(s-1)^{t-1},s)$ for all $i$. Hence the deletion theorem 
implies that 
$\Shi^k_{
+J} \cap H_{\beta_i}^{-k}$ is free with exponents 
$((kh)^\ell)+(d_1,\ldots,d_{\ell-t},(s-1)^{t-1})$. Thus 
$$
|\Shi^k_{
+J}|-
|\Shi^k_{
+J}  \cap H_{\beta_i}^{-k}|
=kh+s-1
$$
for all $i$.

Since $\beta_1+kz,\ldots,\beta_t+kz$ are 
the defining equation of positive roots in $I$ of 
height $kh+s$, it is easy to show that 
they are linearly independent over $\R$. Also, let 
$X:=\cap_{i=1}^t H_{\beta_i}^{-k}$.
If $X \subset H_\alpha^j \in \Shi^k_{+J}$ for some $\alpha \in \Phi^+,\ -k \le j \le k$, then 
puting $z=0$ shows that 
$$
\bigcap_{i=1}^t H_{\beta_i}^0 \subset H_\alpha.
$$
Hence Lemma 4.6 in \cite{ABCHT} implies that $\alpha=\sum_{i=1}^t a_i \beta_i$ with 
$a_i \in \Z_{\ge 0}$. Thus $j=(\sum_{i=1}^t a_i)k$.
Since $-k \le j \le k$, $H_\alpha^j=H_{\beta_i}^{-k}$ for some $i$. 
Hence the three conditions in MAT2 are verified, and we may add hyperplanes 
$H_{\beta_1}^{-k},\ldots,H_{\beta_s}^{-k}$ to 
$\Shi_{+J}^k$, without destroying the freeness, by MAT2, which completes the proof. \owari

Institute of Mathematics for Industry, 
Kyushu University
abe@imi.kyushu-u.ac.jp
%

%
%
%
%

\end{document}